\documentclass[12pt]{gtart}
\usepackage{amsmath}
\usepackage{amssymb}
\usepackage{amsthm}
\usepackage{color}
\usepackage{pinlabel}
\usepackage{subcaption}

\newtheorem{theorem}{Theorem}

\theoremstyle{definition}
\newtheorem{definition}[theorem]{Definition}

\newtheorem{remark}[theorem]{Remark}

\def\C{\mathbb C}

\def\e{\epsilon}

\def\p{\partial}

\begin{document}

\title{Trisections of Lefschetz Pencils}
\authors{ David T. Gay
\footnote{This work was supported by NSF grant DMS-1207721.}}  
\address{Euclid Lab, 160 Milledge Terrace, Athens, GA 30606\\ Department of Mathematics,University
  of Georgia, Athens, GA 30602} \email{d.gay@euclidlab.org}
\begin{abstract}
Donaldson~\cite{DonaldsonLPs} showed that every closed symplectic $4$--manifold can be given the structure of a topological Lefschetz pencil. Gay and Kirby~\cite{GayKirbyTrisections} showed that every closed $4$--manifold has a trisection. In this paper we relate these two structure theorems, showing how to construct a trisection directly from a topological Lefschetz pencil. This trisection is such that each of the three sectors is a regular neighborhood of a regular fiber of the pencil. This is a $4$--dimensional analog of the following trivial $3$-dimensional result: For every open book decomposition of a $3$-manifold $M$, there is a decomposition of $M$ into three handlebodies, each of which is a regular neighborhood of a page.
\end{abstract}

\maketitle

\section{Introduction}

Recall the following two definitions (all manifolds are smooth and oriented and all maps, including coordinate charts, are smooth and orientation preserving):
\begin{definition} \label{D:TLP}
 A {\em topological Lefschetz pencil} on a closed $4$-manifold $X$ is a pair $(B,\pi)$ where $B$ is a nonempty finite collection of points and $\pi: X \setminus B \to S^2$ is a map satisfying:
 \begin{itemize}
  \item \label{I:Base} For each $u \in B$ there is a coordinate chart to $\C^2$ near $u$ and an identification of $S^2$ with $\C P^1$ with respect to which $\pi$ is the standard quotient map $\C^2 \setminus \{0\} \to \C P^1$. Points in $B$ are called {\em base points} and $B$ is the {\em base locus}.
  \item For every critical point $p \in X \setminus B$ of $\pi$, there is a coordinate chart to $\C^2$ near $p$ and a coordinate chart to $\C$ near $\pi(p)$ with respect to which $\pi$ is the map $(z_1,z_2) \mapsto z_1^2 + z_2^2$. These points are called {\em Lefschetz singularities}.
  \item Distinct Lefschetz singularities in $X$ map to distinct values in $S^2$. I.e. each singular fiber has exactly one singularity.
 \end{itemize}
 The {\em genus} $h$ of $(B,\pi)$ is the genus of a (noncompact) regular fiber $\pi^{-1}(q)$, for any regular value $q \in S^2$. If the coordinate charts as in property~\ref{I:Base} above are closed balls $U_1, \ldots, U_b$, then the fiber $\pi^{-1}(q) \cap (X \setminus (\mathring{U}_1 \cup \ldots \cup \mathring{U}_b))$ is called a {\em compact fiber over $q$}, is denoted $F_q$, and when $q$ is regular $F_q$ is a compact surface of genus $h$ with $b$ boundary components.
\end{definition}
Note that the closure of a noncompact regular fiber $\pi^{-1}(q)$ is a smooth closed surface $\pi^{-1}(q) \cup B$ of genus $h$, any two of which intersect transversely and positively at all base points.

\begin{definition} \label{T:Trisection}
 A $(g,k)$--trisection of a closed $4$--manifold $X$ is a decomposition $X=X_1 \cup X_2 \cup X_3$ such that each $X_j \cong \natural^k S^1 \times B^3$, each $X_j \cap X_{j+1} \cong \natural^g S^1 \times B^2$ (indices taken mod $3$) and $X_1 \cap X_2 \cap X_3 \cong \#^g S^1 \times S^1$.
\end{definition}

\begin{theorem} \label{T:FromLPToTrisection}
 Given a $4$--manifold $X$ with a genus $h$ topological Lefschetz pencil $(B,\pi)$, with $b = |B|$ points in the base locus, and $l$ Lefschetz singular fibers, and given any three regular values $y_1,y_2,y_3 \in S^2$, there exists a $(g,k)$ trisection of $X$, $X=X_1 \cup X_2 \cup X_3$, with $g=2h+2b+l-1$ and $k=2h+b-1$, such that each $X_j$ deformation retracts onto a compact fiber $F_{y_j}$ over $y_j$.
\end{theorem}

It will be convenient in describing the relationships between various submanifolds to have the following definition:
\begin{definition}
 We say that $A$ is {\em retraction diffeomorphic} to $B$ when $A$ is a smooth, compact manifold with boundary, $B \subset A$ is a smooth, compact, codimension $0$ submanifold with boundary, and there is a collar structure $(-\e,0] \times \p A$ on a neighborhood of $\p A$ such that $\p B$ is transverse to the $(-\e,0]$ direction, so that $A$ is diffeomorphic to $B$ via a diffeomorphism which is a deformation retraction along the collar.
\end{definition}

\begin{proof}
 Let $u_1, \ldots, u_b \in B$ be the base points, and let $U_1, \ldots, U_b$ be $B^4$ neighborhoods of these points with complex coordinates $(z_1,z_2)$ such that the map $\pi|_{U_i} \co U_i \setminus \{u_i\} \to S^2 = \C P^1$ has the form $(z_1,z_2) \mapsto [z_1:z_2]$, and such that each $U_i$ is a ball of radius $1$. We will first ``trisect'' $Y = X \setminus (\mathring{U}_1 \cup \ldots \cup \mathring{U}_b)$ and then show how to extend this to all of $X$ by filling in with a ``trisection'' of each $U_i$. Until further notice, we now use $\pi$ to refer to $\pi|_Y \co Y \to S^2$.

 Let $p_1, \ldots, p_l \in Y$ be the Lefschetz critical points, with images $q_i = \pi(p_i)$, and let $P_i \ni p_i$, resp. $Q_i \ni q_i$, be $B^4$, resp. $B^2$, neighborhoods such that $\pi(P_i) = Q_i$, with local coordinates on $P_i$ and $Q_i$ with respect to which $\pi$ has the form $(z_1,z_2) \mapsto z_1^2+z_2^2$ and each of $P_i$ and $Q_i$ have radius $1$. Assume that none of the $Q_i$'s contain any of the given regular values $y_1$, $y_2$ or $y_3$, or the north or south poles of $S^2$. Let $y_+$ refer to the north pole of $S^2$ and let $y_-$ refer to the south pole.

 Split $S^2$ into three bigons $S^2 = A_1 \cup A_2 \cup A_3$ intersecting at the north and south poles as in Figure~\ref{F:TrisectS2}, with $y_j \in \mathring{A}_j$, with each $q_i \in A_1 \cap A_3$, each $Q_i$ contained in $A_1 \cup A_3$, and $A_1 \cap Q_i = \{z \mid \mathop{Im}(z) \leq 0\}$ and $A_3 \cap Q_i = \{z \mid \mathop{Im}(z) \geq 0\}$. 
 \begin{figure}[h]
 \labellist
 \tiny\hair 2pt
 \pinlabel $y_+$ [bl] at 81 118
 \pinlabel $q_1$ [br] at 67 96
 \pinlabel $Q_1$ at 73 93
 \pinlabel $q_2$ [r] at 51 71
 \pinlabel $Q_2$ at 58 66
 \pinlabel $A_1$ at 94 34
 \pinlabel $A_2$ at 105 119
 \pinlabel $A_3$ at 28 89
 \endlabellist
 \centering
 \includegraphics[width=2in]{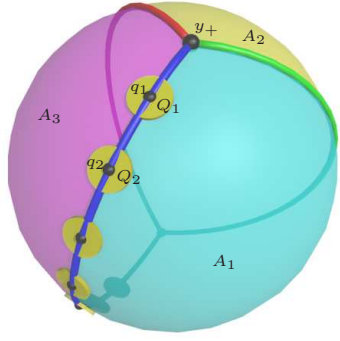}
 \caption{\label{F:TrisectS2} ``Trisection'' of $S^2$.}
 \end{figure}

 Let $P_i' = \{|(z_1,z_2)| \leq 1/2\} \subset P_i$ be the $4$--ball of radius $1/2$ inside the $4$--ball $P_i$ of radius $1$. We carefully connect each $P_i'$ to $\pi^{-1}(A_2)$ as follows. For each $i = 1, \ldots, l$, consider an embedding of a rectangle $\rho_i \co [-\e,\e] \times [\e,1] \to S^2$ (for suitably small $\e > 0$) as illustrated in Figure~\ref{F:RectangleInS2}, such that:
 \begin{enumerate}
  \item For all $t \in [\e,1/2+\e]$, $\rho_i(s,t) = s + t i \in Q_i$, using the complex coordinate on $Q_i$ with respect to which $\pi$ is $(z_1,z_2) \mapsto z_1^2 + z_2^2$.
  \item $\rho_i(0,3/4) = y_+$, the north pole
  \item For $t \in [\e,3/4]$, $\rho_i(0,t)$ is an embedded path in $A_1 \cap A_3$
  \item For $t \in [3/4,1]$, $\rho_i(0,t)$ is embedded in $A_2$.
  \item The decomposition $S^2 = A_1 \cup A_2 \cup A_3$ splits the rectangle $R = [-\e,\e] \times [0,1]$ into two trapezoids $R_1$ and $R_3$ and a pentagon $R_2$, with $R_1 \cap R_3 = \{0\} \times [0,3/4]$, $R_1 \subset [-\e,0] \times [0,1]$ and $R_3 \subset [0,\e] \times [0,1]$, exactly as in Figure~\ref{F:RectangleInS2}.
 \end{enumerate}
 \begin{figure}[h]
 \labellist
 \tiny\hair 2pt
 \endlabellist
 \centering
 \includegraphics[width=1.5in]{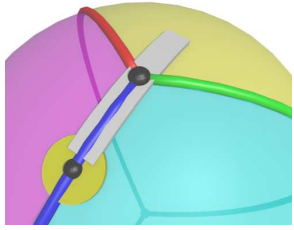}
 \caption{\label{F:RectangleInS2} The rectangle $\rho_1$ in $S^2$.}
 \end{figure}
 Now lift this to an embedding $\tau_i \co B^2_\e \times [-\e,\e] \times [1/4,1] \to Y$ such that $\pi(\tau_i(x,y,s,t)) = \rho_i(s,t)$, arranging also that all the $\tau_i$'s have disjoint images. Finally, let $N_i = \tau_i(B^2_\e \times [-\e,\e] \times [1/4,1]) \cup U_i$; this is a thin finger sticking out of $\pi^{-1}(X_2)$ with the ball $U_i$ stuck to the end. The part of each finger that is outside $\pi^{-1}(A_2)$ is evenly split, with half in $\pi^{-1}(A_1)$ and half in $\pi^{-1}(A_3)$.

 Now we define the trisection of $Y$ as follows:\begin{enumerate}
  \item $Y_1 = \pi^{-1}(A_1) \setminus (\mathring{N}_1 \cup \ldots \cup \mathring{N}_l)$
  \item $Y_2 = \pi^{-1}(A_2) \setminus (\mathring{N}_1 \cup \ldots \cup \mathring{N}_l)$
  \item $Y_3 = \pi^{-1}(A_3) \cup N_1 \cup \ldots \cup N_l$
 \end{enumerate}

 It is not hard to see that each $Y_j$ is retraction diffeomorphic to a tubular neighborhood of $\pi^{-1}(y_j)$, and is hence diffeomorphic to $B^2 \times F_{h,b}$, where $F_{h,b}$ is a compact surface of genus $h$ with $b$ boundary components. Note that $B^2 \times F_{h,b} \cong \natural^{2h+b-1} S^1 \times B^3$. We need to understand the pairwise and triple intersections.

 The pairwise intersection $Y_1 \cap Y_2$ is the boundary connected sum of $\pi^{-1}(A_1 \cap A_2) \cong F_{h,b} \times [0,1]$ with one solid torus in the boundary of each $P_i'$. The Morse function $\mathop{Im}(z_1^2+z_2^2)$ on $\p P_i'$ splits $\p P_1' \cong S^3$ into two solid tori, one in $A_1$ and one in $A_3$, and the one in $A_1$ is connected back to $\pi^{-1}(A_1 \cap A_2)$ along the boundary of the tubular neighborhood $\nu_i$. Thus $Y_1 \cap Y_2 \cong \natural^{2h+l+b-1} S^1 \times B^2$, and the same argument holds for $Y_3 \cap Y_2$.
 
 To understand $Y_1 \cap Y_3$, consider a meridian line $m \subset S^2$ from the north pole $y_+$ to the south pole $y_-$, parallel and close to $A_1 \cap A_3$, just inside $A_1$ and outside the $Q_i$'s. Then $\pi^{-1}(m) \cong \pi^{-1}(y_+) \times [0,1] \cong F_{h,b} \times [0,1] \cong \natural^{2h+b-1} S^1 \times B^2$; this is illustrated in Figure~\ref{F:BaseAndBaseXI}. 
 \begin{figure}[h]
    \centering
    \begin{subfigure}[b]{0.4\textwidth}
        \centering
        \includegraphics[width=2in]{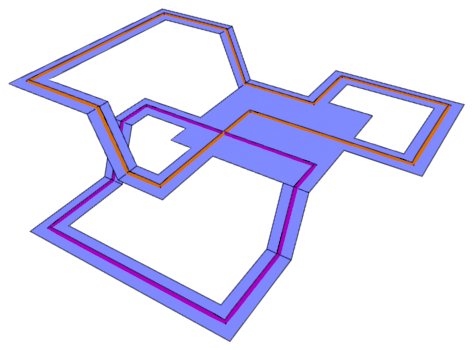}
    \end{subfigure}
    \begin{subfigure}[b]{0.4\textwidth}
        \centering
        \includegraphics[width=2in]{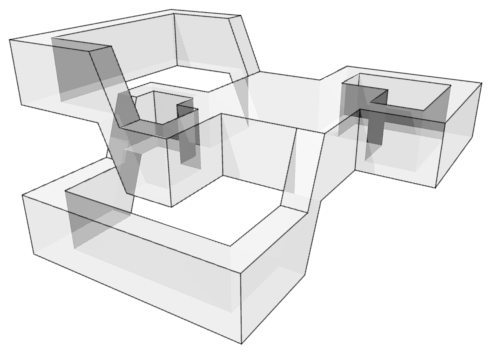}
    \end{subfigure}
 \caption{\label{F:BaseAndBaseXI} On the left, $\pi^{-1}(y_+) \cong F_{h,b}$; here we have also drawn two curves which will be important in the following figure. On the right, $\pi^{-1}(m) \cong F_{h,b} \times [0,1]$.}
 \end{figure}
 If we rotate $m$ towards the meridian line $A_1 \cap A_3$, we see $Y_1 \cap Y_3$ as diffeomorphic to the complement in $\pi^{-1}(m) \cong F_{h,b} \times [0,1] \cong \natural^{2h+b-1} S^1 \times B^2$ of a collection of $l$ solid tori, each connected to $\pi^{-1}(y_+) \cong F_{h,b} \times \{0\}$ by a tube. Each of these solid tori in $F_{h,b} \times [0,1]$ is a tubular neighborhood of a simple closed curve in a surface $F_{h,b} \times \{t_i\}$ for some $t_i \in [0,1]$ (corresponding to $q_i \in A_1 \cap A_3$); these solid tori (horizontal) and connecting tubes (vertical posts) are illustrated in Figure~\ref{F:BaseXIWithCurvesAndPosts}.
 \begin{figure}[h]
 \labellist
 \tiny\hair 2pt
 \endlabellist
 \centering
 \includegraphics[width=2in]{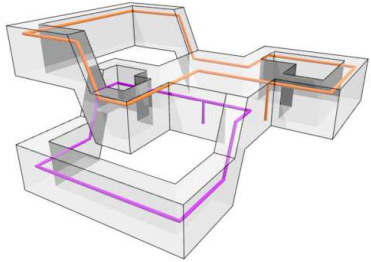}
 \caption{\label{F:BaseXIWithCurvesAndPosts} Solid tori and posts in $F_{h,b} \times [0,1]$; the complement of the solid tori and posts is $Y_1 \cap Y_3$.}
 \end{figure}
 The complement of these solid-tori-with-posts is in fact a solid handlebody diffeomorphic to $\natural^{2h+l+b-1} S^1 \times B^2$ for the following reason: Removing each solid-torus-with-post is the same as removing a single thickened arc which starts on $F_{h,b} \times \{0\}$, goes up to $F_{h,b} \times \{t_i\}$ parallel to the post, goes around the loop in $F_{h,b}$ without quite closing up, then goes back to $F_{h,b} \times \{0\}$ parallel to the tube. Removing these thickened arcs one by one, starting from the highest one, the largest value of $t_i$, corresponds to adding $1$--handles one by one to $F_{h,b} \times [0,1]$ because each of these arcs is boundary parallel in the complement of the higher arcs; see Figure~\ref{F:WhyItsAHandlebody}
 \begin{figure}[h]
    \centering
    \begin{subfigure}[b]{0.4\textwidth}
        \centering
        \includegraphics[width=2in]{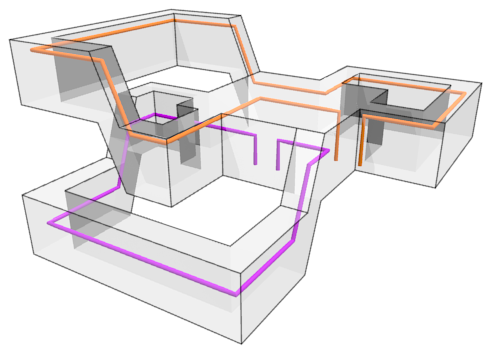}
    \end{subfigure}
    \begin{subfigure}[b]{0.4\textwidth}
        \centering
        \includegraphics[width=2in]{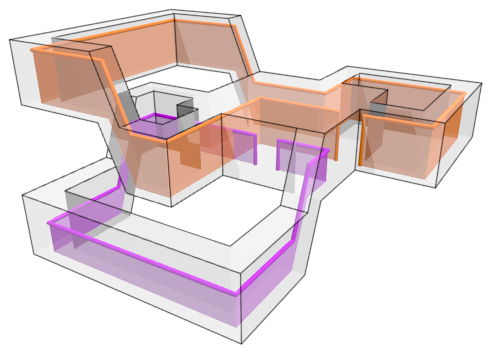}
    \end{subfigure}
 \caption{\label{F:WhyItsAHandlebody} Seeing that $Y_1 \cap Y_3$ is a handlebody. On the left we see $Y_1 \cap Y_3$ as the complement of arcs in $F_{g,b} \times [0,1]$, and on the right we see that these arcs can be ordered so that each is boundary parallel in the complement of the preceding arcs.}
 \end{figure}
 
 The triple intersection $Y_1 \cap Y_2 \cap Y_3$ lies entirely in $\pi^{-1}(A_1 \cap A_3)$, and is equal to the connected sum of $\pi^{-1}(A_1 \cap A_2 \cap A_3) \cong F_{h,b} \times \{y_+,y_-\}$ with $l$ copies of $T^2$ (the connected sum happening on the $F_{h,b} \times \{y_+\}$ component). Each of the $T^2$'s is a separating $T^2$ in $\p P_i \cong S^3$. Thus $Y_1 \cap Y_2 \cap Y_3 \cong F_{h,b} \amalg F_{h+l,b}$.

 Now we need to fill in the neighborhoods $U_i$ of the base points $u_i \in B$. We will do the same thing in each $U_i$, so we will drop the subscript $i$ and simply look at one neighborhood $U = B^4$. On $\p B^4 = S^3$, we have $S^3 = (S^3 \cap Y_1) \cup (S^3 \cap Y_2) \cup (S^3 \cap Y_3)$ where each $S^3_j = S^3 \cap Y_j$ is $(\pi|_{S^3})^{-1}(A_j)$. But $\pi|_{S^3}$ is just the Hopf fibration $S^3 \to S^2$, and we see a decomposition of $S^3$ into three solid torus neighborhoods, $S^3_1$, $S^3_2$ and $S^3_3$, of three Hopf fibers, and these three neighborhoods have triple intersection equal to a Hopf link and pairwise intersections equal to annulus pages of the Hopf open book decomposition of $S^3$ for this Hopf link. This is illustrated in Figure~\ref{F:S3Trisection}. We now need to extend this decomposition to all of $B^4$.
 \begin{figure}[h]
    \centering
    \begin{subfigure}[b]{0.4\textwidth}
        \centering
        \includegraphics[width=1in]{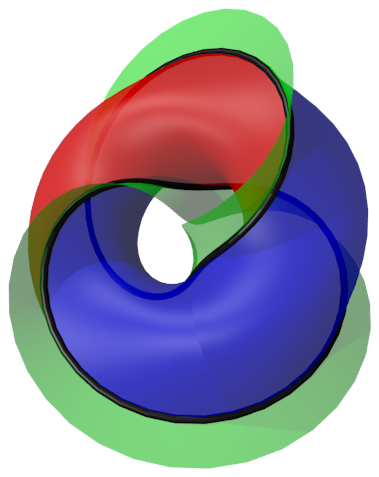}
       \subcaption{Trisection of $S^3 = \p U_i$}
        \label{F:S3Trisection}
    \end{subfigure}
    \begin{subfigure}[b]{0.4\textwidth}
        \centering
        \includegraphics[width=1in]{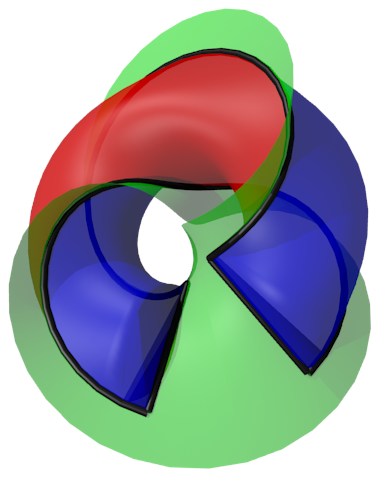}
       \subcaption{After pinching.}
        \label{F:PinchSolidTorus}
    \end{subfigure}
 \caption{The trisection of $S^3 = \p U_i$ into three solid tori. Here we focus on what is happening near one of these solid tori. The colors identify the three annular pages of the Hopf open book which are also the three pairwise intersections $S^3_1 \cap S^3_2$, $S^3_2 \cap S^3_3$ and $S^3_3 \cap S^3_1$. Black is the triple intersection, a Hopf link, and the solid tori themselves are the regions between the annuli. On the right we have pinched one of the solid tori $S^3_j$ to a ball and a rectangle; here red and blue have changed from annuli to disks while green has changed from an annulus to a genus $1$ surface with $1$ boundary component}
 \end{figure}

 We describe the extension of $S^3 = S^3_1 \cup S^3_2 \cup S^3_3$ to $B^4 = B^4_1 \cup B^4_2 \cup B^4_3$ as a $1$--parameter family of decompositions of concentric $S^3$'s in $B^4$, i.e. as a movie as we shrink the $S^3$ from $\p B^4$ to the center. As we shrink the $S^3$, we will consecutively pinch each solid torus $S^3_j \cong S^1 \times B^2$ to a  $3$--ball with a rectangle attached, as in Figure~\ref{F:PinchSolidTorus}. The $3$--ball is what we see of $B^4_j$ at this radius, while the rectangle becomes part of $B^4_{j-1} \cap B^4_{j+1}$ at this radius (and further inwards). After pinching each solid torus, we are left with a radius at which we have three $B^3$'s intersecting in pairs along $B^2$'s and with triple intersection an unknot, and this continues to the center. The upshot is that each $B^4_j$ is a copy of $S^1 \times B^3$ with half of its boundary ($S^1 \times B^2$) on $\p B^4 = S^3$ and the other half in the interior. The pairwise intersections $B^4_j \cap B^4_{j+1}$ are genus $2$ handlebodies with their boundaries split as an annulus $S^3_j \cap S^3_{j+1}$ on $\p B^4$ and a genus $1$ surface with $2$ boundary components properly embedded in $B^4$. The triple intersection is this properly embedded genus $1$ surface with $2$ boundary components.

 So now we extend the decomposition $Y=Y_1 \cup Y_2 \cup Y_3$ to an honest trisection $X=X_1 \cup X_2 \cup X_3$ by attaching each $B^4_j$ in each $B^4 = U_i$ to $Y_j$. This does not change the topology of the $4$--dimensional pieces, and in fact we see that $X_j$ is retraction diffeomorphic to $Y_j$ via a retraction along each $B^4_j$, and thus $X_j \cong \natural^{2h+b-1} S^1 \times B^3$. However, $X_j \cap X_{j+1}$ is not diffeomorphic to $Y_j \cap Y_{j+1}$; the difference is that we have attached genus $2$ handlebodies along annuli in $\p (Y_j \cap Y_{j+1})$, which really means that we have attached one $3$--dimensional $1$--handle for each base point. Thus $X_j \cap X_{j+1} \cong \natural^{2h+2l+b-1} S^1 \times B^2$, and similarly $X_1 \cap X_2 \cap X_3 \cong \#^{2h+2l+b-1} S^1 \times S^1$. (Alternatively, for $X_1 \cap X_2 \cap X_3$, recall that $Y_1 \cap Y_2 \cap Y_3 \cong F_{h,b} \amalg F_{h+l,b}$. Filling in the trisection of one $U_i$ connects one boundary component of $F_{h,b}$ with one boundary component of $F_{h+l,b}$ via a genus $1$ surface with $2$ boundary components, ultimately producing a genus $h+(h+l)+l+(b-1)$ surface.)
 
 Finally, note that each $Y_j$ deformation retracts onto the compact fiber $(\pi|_Y)^{-1}(y_j) \subset Y$, and hence so does each $X_j$. But as constructed, $X_j$ is not a tubular neighborhood of the full noncompact fiber $\pi^{-1}(y_j) \subset X$.

\end{proof}

\begin{remark}
 This theorem as stated is not in fact sensitive to the orientations involved in the local models in the definition of ``Lefschetz pencil'', so that the theorems holds for ``achiral Lefschetz pencils'' as well. The achirality may be interpreted as removing the orientation constraints at either the base points or the critical points.
\end{remark}

\begin{remark}
 We have stated this theorem in the honest Lefschetz pencil case, not the achiral case, simply because the main motivation is to work toward a notion of a ``symplectic trisection''.
\end{remark}

\begin{remark}
 Running the construction above in the case of the pencil of lines in $\C P^2$ yields the standard trisection of $\C P^2$ into three $B^4$'s meeting in pairs along solid tori, coming from the toric structure on $\C P^2$, as described in~\cite{GayKirbyTrisections}.
\end{remark}

\begin{remark}
 This theorem does {\em not} hold for Lefschetz fibrations. The essential problem is that we have closed fibers, and a regular neighborhood of a closed surface is not a $1$--handlebody.
\end{remark}

\begin{remark}
 We have not addressed the question of how to draw a trisection {\em diagram} for a trisection coming from a Lefschetz pencil. Also note that the phrase ``vanishing cycle'' does not appear in the proof above. Obviously a description of the pencil in terms of vanishing cycles will be sufficient to reproduce the trisection diagram. This remains to be worked out in detail, but a sketch of the method is as follows: In fact, vanishing cycles did appear in this proof, they are the simple closed curves in parallel copies $F_{h,b} \times \{t_i\} \subset F_{h,b} \times [0,1]$ of the compact fiber $F_{h,b}$, as in Figure~\ref{F:BaseXIWithCurvesAndPosts}. Each vanishing cycle is connected to $F_{h,b} \times \{0\}$ with a post, and the complement is $Y_1 \cap Y_3$. We can then see directly the surface $Y_1 \cap Y_2 \cap Y_3$ as part of the boundary of this complement (the top by itself, and the bottom with the boundaries of the posts and solid tori). From the picture of $Y_1 \cap Y_3$ we can see which curves bound disks in $Y_1 \cap Y_3$. To see which curves bound disks in $Y_1 \cap Y_2$ and $Y_2 \cap Y_3$, we have meridian and $(1,1)$--curve pairs on each of the solid tori, and the remaining curves come in parallel pairs. Finally, when we extend to the trisection $X=X_1 \cup X_2 \cup X_3$, we get meridian, longitude and $(1,1)$-curve on each of the new $S^1 \times S^1$ summands of the central surface.
\end{remark}

\bibliographystyle{plain}
\bibliography{TrisectionsOfLPs}

\end{document}